\documentclass{jssc}
\usepackage{color}
\def\textsubscript#1%
{$_{\text{#1}}$}

\def\cdd{\mbox{\boldmath$\cdot$}~}
\def\dfrac{\displaystyle\frac}
\input{amssym.def}
\include{graphix}

\makeatletter
\def\@oddfoot{\hfill}
\newcount\shumeicount
\def\setshumei#1#2#3{%
  \shumeicount=\count0
  \def\@oddhead{%
    \raise-5pt\hbox to0pt{\vrule width\hsize height 0pt depth 0.4pt\hss}\relax
    \ifnum \shumeicount=\count0
      \raise-7pt\hbox to0pt{\vrule width\hsize height 0pt depth 0.4pt\hss}\relax
      #1
    \else
      \ifodd\count0
        #2
      \else
        #3
       \fi
     \fi
  }%
}
\makeatother
\makeatletter
\def\@oddfoot{\hfill}
\newcount\shujiaocount
\def\setshujiao{%
  \shujiaocount=\count0
  \def\@oddfoot{%
      \ifodd\count0
      \else
      \fi
  }%
}
\makeatother
\def\title#1#2#3#4{{
  \vspace*{0.3cm}
  \begin{flushleft} \Large\bf #1\end{flushleft}
  \vspace*{-0.2cm}
      \begin{flushleft}
      \bf #2
      \end{flushleft}
      \footnotetext{\hspace{-6mm} #3\\ #4}}}
\def\dshm#1#2#3#4
{\setshumei{J Syst Sci Complex (#1) #2:
{\thepage--\pageref{LastPage}}\hfill}
            {\hfill {\small #3}\hfill\hbox to0pt{\hss\thepage}}
            {\hbox to0pt{\thepage\hss}\hfill {\small #4}\hfill
            }
            \setshujiao}
\def\drd#1#2
{{\vskip 1cm\small \begin{flushleft}
 #1 \\
 #2\\
\copyright The Editorial Office of JSSC \&  Springer-Verlag GmbH
Germany 2018
\end{flushleft}}}
\def\tilde{\widetilde}
\def\bar{\overline}
\def\epsilon{\varepsilon}

\begin{document}

%*************************************************************************************************************
% \biaoti{THE CAPITALIZED TITLE OF YOUR ARTICLE$^*$}{The list of authors' names with the LAST NAME capitalized
% and the authors' names should be separated by "\cdd"}{the first author's name \\ the first author's affiliation
% and Email address\\ the second author's name\\ the second author's affiliation. More can be listed like this.}
% {$^*$ The titles and numbers of the foundations that support this article.}
%*************************************************************************************************************

%Yuke Shi$^{1,2}$, Wei Zhang$^{1}$, Aiyi Liu$^{3}$, Qizhai Li$^{1,2,*}$
%
%$^1$LSC, NCMIS, Academy of Mathematics and Systems Science, Chinese Academy of Sciences, Beijing 100190, China.
%
%$^2$University of Chinese Academy of Sciences, Beijing 100049, China.
%
%$^3$Biostatistics and Bioinformatics Branch, {\it Eunice Kennedy Shriver} National Institute of Child Health and Human Development, National Institutes of Health, Bethesda, MD 20847, USA.
%
%$^*$ Corresponding author. E-mail: liqz@amss.ac.cn

\title{Distance-based regression analysis for measuring associations$^*$}%%%   Main Title of your paper  %%%
{\uppercase{Shi} Yuke  \cdd \uppercase{Zhang}
Wei \cdd \uppercase{Liu } Aiyi \cdd \uppercase{Li} Qizhai }%%% The names of the authors  %%%
{\uppercase{Shi} Yuke \\
 \emph{LSC, Academy of Mathematics and Systems Science, Chinese Academy of Sciences, Beijing 100190, China; University of Chinese Academy of Sciences, Beijing 100049, China. } Email:shiyuke@amss.ac.cn \\   % Academy of Mathematics and Systems Science, Chinese Academy of Sciences, Beijing $100190$, China
 \uppercase{Zhang} Wei \\
 \emph{LSC, Academy of Mathematics and Systems Science, Chinese Academy of Sciences, Beijing 100190, China.}  Email:zhangwei@amss.ac.cn \\
\uppercase{Liu } Aiyi \\
 \emph{Biostatistics and Bioinformatics Branch, {\it Eunice Kennedy Shriver} National Institute of Child Health; Human Development, National Institutes of Health, Bethesda, MD 20847, USA.}
 Email:liua@mail.nih.gov \\
\uppercase{Li} Qizhai (Corresponding author) \\
 \emph{LSC, Academy of Mathematics and Systems Science, Chinese Academy of Sciences, Beijing 100190, China; University of Chinese Academy of Sciences, Beijing 100049, China.}  Email:liqz@amss.ac.cn \\
} %%% The address of the authors  %%%
{$^*$This work was partially supported by Beijing Natural
	Science Foundation (Z180006).\\
%{$^\diamond${\it This paper was recommended for publication by Editor . }}
}

%*************************************************************************************************************
%The submission date of your article. For example: \drd{Received: June 8, 2006}
%*************************************************************************************************************
\drd{DOI: }{Received: x x 20xx}{ / Revised: x x 20xx}

%*************************************************************************************************************
% The page header of the article.
% \dshm{Year}{Volume}{The capitalized RUNNING HEAD of your article with less than 48 letters}{The capitalized
% AUTHORS list with $\cdot$ separating different names or one can type "The name of the first author et al."
% if there are more than 4 authors.}
%*************************************************************************************************************

%%%%\dshm{20XX}{XX}{A TEMPLATE FOR JOURNAL}{\uppercase{Surname1 Firstname1} $\cdd$ \uppercase{Surname2 Firstname2} $\cdd$\uppercase{Surname3 Firstname3}}

%*************************************************************************************************************
% \dab{The abstract}{Keywords}
%*************************************************************************************************************
%-------------------------------------------------------------------------
\Abstract{Distance-based regression model, as a nonparametric multivariate method, has been widely used to detect the association between
	variations in a distance or dissimilarity matrix for outcomes and predictor variables of interest in genetic association studies, genomic analyses, and many other research areas. Based on it, a pseudo-$F$ statistic which partitions the variation in distance  matrices is often constructed to achieve the aim. To the best of our knowledge, the statistical properties of the pseudo-$F$ statistic has not yet been well established in the literature. To fill this gap, we study the asymptotic null distribution of the pseudo-$F$ statistic and show that it is asymptotically equivalent to a mixture of chi-squared random variables. Given that the pseudo-$F$ test statistic has unsatisfactory power when the correlations of the response variables are large, we propose a square-root $F$-type test statistic which replaces the similarity matrix with its square root. The asymptotic null distribution of the new test statistic and power of both tests are also investigated.
	Simulation studies are conducted to validate the asymptotic distributions of the tests and demonstrate that the proposed test has more robust power than the pseudo-$F$ test. Both test statistics are exemplified with a gene expression dataset for a prostate cancer pathway.
}      % the abstract

\Keywords{Asymptotic distribution, Chi-squared-type mixture,  Nonparametric test, Pseudo-$F$ test, Similarity matrix.}        % the keywords

%\MRSubClass{05B05, 05B25, 20B25}      % MR(2000) Subject Classification

%\baselineskip 15pt

\section{Introduction}
Distance-based regression model, as a multivariate method testing  the relationship between variation in a distance or dissimilarity matrix of outcomes and predictor variables of interest,
has been widely used in a variety of applications including  genetics \citep{Han2010}, genomics \citep{Zapala2006, Nievergelt2007}, 
microbiome \citep{Wu2011, Norman2015, Liang2015, Wang2019}, and other research area such as geoscience \citep{Molari2018},environmental science \citep{White2020}, and oceanography \citep{Bertocci2012, Consoli2013}.
Unlike conventional multivariate methods, the starting point of distance-based regression model  is pairwise distance or dissimilarity of subjects instead of individuals' observations. This model was originally introduced by \citet{McArdle2001} to analyze simultaneous responses of multiple species to several factors in ecological experimental designs. Subsequently, with the development of high-throughput technologies, it has been received considerable attentions in high-dimensional data analysis. For example, \citet{Wessel2006} applied distance-based regression model to relate variations measured by a genomic distance  for multiple phenotypes with multiple variants; \citet{Zapala2006} proposed a multivariate regression approach based on a distance matrix to detect associations between expression patterns on groups of genes and related  predictor variables; \citet{Chen2012} developed a distance-based statistical test based on generalized UniFrac distances to detect the association of microbiome composition and environment covariates. \citet{Gambi2020} applied it to detect the impaction of historical sulfide mine tailings discharge on meiofaunal assemblages.

By partitioning the variation in a distance or dissimilarity matrix, a pseudo-$F$ test statistic, analogous to the conventional univariate $F$ statistic  in ANOVA or  general linear models,  is constructed to detect the association between multiple response and predictor variables \citep{Reiss2010}.
However, to the best of our knowledge, the distribution  or asymptotic distribution of the pseudo-$F$ test statistic is unknown with only some results under normal assumption and Euclidean or Mahalanobis distance measures \citep{LiJL2019}.
For the general cases, it has not been studied yet in the literature.
Resampling procedures are often used to calculate the statistical significance  of the pseudo-$F$ test statistic \citep{McArdle2001, Wessel2006, Zapala2006, Li2009, Han2010, Chen2012}. However, they are generally computationally intensive, especially for high-dimensional data.

To address this problem, we study the asymptotic distribution of the pseudo-$F$ statistic. We show that under the null hypothesis, the pseudo-$F$ statistic is asymptotically equivalent to  a mixture of Chi-squared random variables, with the weights  proportional to the eigenvalues of the similarity matrix. Such Chi-squared-type mixture variables could have a large critical value when the correlation coefficients of the response variables are large, since in this case the eigenvalues of the similarity matrix vary widely. Thus the pseudo-$F$ statistic may suffer from a substantial power loss. So we propose a square-root $F$-type test statistic, which
reduces the difference in the eigenvalues of the similarity matrix by replacing the similarity matrix with its square root. The asymptotic null distribution is also established for the new test statistic.

The rest of the article is organized as follows. We introduce the pseudo-$F$ test statistic and the proposed test statistic in Section 2. Main theoretical results for the test statistics are given in Section 3. Simulation studies are conducted to validate the asymptotic distribution of the pseudo-$F$ and square-root $F$-type test  statistics in Section 4. In Section 5, we illustrate the two tests  with gene expression data from a prostate cancer pathway study. Some conclusions are drawn in Section 6. Technical details are relegated to the Appendix.

\section{Test statistics}

We introduce the pseudo-$F$ test statistic by taking a similarity matrix as the starting point since the distance matrix can be easily transferred to a similarity matrix.
Let $S=(s_{ij})_{n\times n}$ be a similarity matrix for $k$ response variables among $n$ subjects, where $s_{ij}$ is the pairwise similarity between the subjects $i$ and $j$ constructed  based on  a positive definite kernel $\psi(y_i,y_j)$, and $y_i$ and $y_j$ are two high-dimensional independent observations of the responses variable $\mathbb Y$, $i,j=1,\cdots, n$. Let $\mathbb X=(\mathbb X_1,\cdots,\mathbb X_m)^\top$ be $m$ predictor variables of interest. Denote the observation of $\mathbb X$ for the $i$th subject by $x_i=(x_{i1},\cdots,x_{im})^\top$, $i=1,\cdots, n$. Write $X= (x_1,\cdots,x_n)^\top$. We refer to the regression model relating $S$ and $X$ as a distance-based regression model, denoted by $S\sim X$. Our interest is in testing the null hypothesis $H_0$: there is  no association between $\mathbb Y$ and $\mathbb X$.
A $F$-type test statistic proposed by \citet[page.~3]{McArdle2001}  can  be  used  to do it, which is given by
\begin{equation*}
T_{pseudo}=\frac{\text{tr}\big(H_{X}{{\bf H}S{\bf H}}/m \big)}{\text{tr}\big(({\bf I}_n - {H_X}){\bf H}S{\bf H}/(n-m)\big)},
\end{equation*}
where $H_X= X( X^\top X)^{-1}X^\top$ is the traditional hat" matrix, ${\bf H}={\bf I}_n-{\bf 1}_n{\bf 1}_n^\top/n$,  ${\bf I}_n$ is the $n\times n$ identity matrix, ${\bf 1}_n$ is an $n$-dimensional column vector with all elements being 1, and $\text{tr}(\cdot)$ represents the trace of a matrix. Since the distribution or asymptotic distribution of $T_{pseudo}$ is unknown, its statistical significance is conventionally calculated using resampling procedures, which is usually computationally intensive, especially when $n$ is large. One of the goals of this article is to establish the asymptotic distribution for the pseudo-$F$ test statistic $T_{pseudo}$.

The numerical results in Section 4 below  show that $T_{pseudo}$ is sensitive to the correlation matrix of the response variables. To boost the power when the response variables are moderately or highly correlated, we propose a square-root $F$-type test statistic, which uses the square root of the similarity matrix. The square-root $F$-type statistic is defined as
$$T_{sqrt}=\frac{\text{tr}\big(H_X({\bf H}S{\bf H})^{1/2}/m \big)}{\text{tr}\big(({\bf I}_n - H_X)({\bf H}S{\bf H})^{1/2}/(n-m)\big)},$$
where the square root $B$ of a matrix $A$ means $A=B\times B$. Consider the  eigenvalues $\{\lambda_i\}_{i=1}^n$ and eigenfunctions $\{v_i\}_{i=1}^n$ of $\tilde{S}={\bf H}S{\bf H}$ and    $\tilde{S}^{1/2}=(v_1,\cdots,v_n)diag(\lambda_{1}^{1/2},\cdots,\lambda_{n}^{1/2})(v_1,\cdots,v_n)^\top$. 
The numerator of $T_{pseudo}$ and $T_{sqrt}$ now can be intuitively expressed as $\sum_{i=1}^n\lambda_i v_i^\top \tilde{v}_i/m$ and  $\sum_{i=1}^n\lambda_i^{1/2}  v_i^\top \tilde{v}_i/m$ where the eigenvalues of $H_x$ are all $1$ and  the eigenfunctions of $H_x$ are $\{\tilde{v}_i\}_{i=1}^n$. Superficially, when the response variables are moderately correlated,  our square-root method may increase the weight of valuable factors to boost the power.

\section{Theoretic properties and implementation}\label{MR}
In this section, we derive the asymptotic distributions of the pseudo-$F$  and the proposed  test statistics.   Assume that $E(\mathbb X) =\mu$ and $\text{cov}(\mathbb X)=\Delta=(\delta_{ij})_{m\times m}$. It thus follows that $E(\mathbb X \mathbb X^\top)=\Delta + \mu \mu^\top\triangleq\tilde\Delta=(\tilde\delta_{ij})_{m \times m}$.   Denote the eigenvalues of $\tilde{S}={\bf H}S{\bf H} $  by  $\lambda_{1}\geq  \lambda_{2}\geq  \cdots \geq  \lambda_{n} \geq  0$. Besides, let $\lambda^*_1\geq   \cdots \geq  \lambda^*_n \geq  0$ be the solutions to  the  equation
$\int \psi_0 (y_1,y_2)  u_i (y_1)p(y_1)dy_1=\lambda^*_i  u_i(y_2),$
where $\psi_0(y_1,y_2)=\psi(y_1,y_2)-E_{y_1}(\psi(y_1,y_2))- E_{y_2}(\psi(y_1,y_2))+E_{y_1,y_2}(\psi(y_1,y_2))$ and $ u_i(y_1)$ is the eigenfunction of the kernel $\psi_0(y_1,y_2)$ corresponding to $\lambda^*_i$ for $i=1,\cdots,n$.
Throughout, we assume that there exists positive values $c_0$ and $c_1$,  such that $E(\mathbb X_j^4) \leq  c_0$ and $E(s_{ii}^2)\leq  c_1$ for $j=1,\cdots,m$ and $i=1,\cdots,n$.  The symbol $``\stackrel{d}\longrightarrow"$ represents ``converges in distribution to" and $``\stackrel{p}\longrightarrow"$ means ``converges in probability to".

\subsection{Asymptotical null distributions}

In the following two lemmas, we give the asymptotic null distribution of the numerator  and denominator of the pseudo-$F$ test statistic, respectively.

\begin{lemma}\label{lem1}
	Under the null hypothesis $H_0$, the numerator of $T_{pseudo}$, $\text{tr}\big(H_X{\bf H}S{\bf H}\big)$, has the same asymptotic distribution as $\frac{1}{n}\sum_{i=1}^{n} \lambda_i \xi_{i}$ when $n\rightarrow\infty$, where $\xi_i=\varpi_i+\varrho_i/(1+\mu^\top \Delta^{-1}\mu) $, $\varpi_1,\ldots,\varpi_n$ are independently identical distributed (\textit{iid}) chi-squared random variables  with  $m-1$ degrees of freedom, $\varrho_1, \cdots, \varrho_n$ are \textit{iid} chi-squared random variables  with 1 degrees of freedom, i.e.,  $\varpi_i\stackrel{\text{iid} }\thicksim \chi^2_{m-1}$ and $\varrho_i\stackrel{\text{iid} }\thicksim \chi^2_1$, and $\varpi_i$ and $\varrho_i$ are independent, $i=1,\cdots, n$. In particular, when $\mu= \bm 0_m$ an $m$-dimensional column vector with all 0 units, $\xi_i\stackrel{\text{iid} }\thicksim \chi^2_{m}$.
\end{lemma}

\begin{lemma}\label{lem2}
	As $n\rightarrow\infty$, the denominator of $T_{pseudo}$,
	\begin{align*}
	\text{tr}\big(({\bf I}_n - H_X){\bf H}S{\bf H}/(n-m)\big) \stackrel{p} \longrightarrow E(s_{11})-E(s_{12}),
	\end{align*}
	where $E(s_{11})= E(\psi(y_1, y_1))$, $E(s_{12})= E(\psi(y_1, y_2))$, and $ y_1$ and $ y_2$ are two independent observations of $\mathbb Y$.
\end{lemma}

All those proofs are presented in Appendix. Both lemmas lead to the following theorem concerning the asymptotic distribution of  the pseudo-$F$ test statistic.

\begin{theorem}\label{thm1}
	Under the null hypothesis $H_0$,  $T_{pseudo}$ has the same asymptotic distribution with $m^{-1}\sum_{i=1}^{n} w_i \xi_{i}$ when $n\rightarrow\infty$, where $\xi_i=\varpi_i+\varrho_i/(1+\mu^\top \Delta^{-1} \mu) $, $\varpi_i\stackrel{\text{iid} }\thicksim \chi^2_{m-1}$ and $\varrho_i\stackrel{\text{iid} }\thicksim \chi^2_1$
	are mutually independent, $w_i=\lambda_i/\sum_{j=1}^{n} \lambda_j$, and $\{\lambda_1,\cdots, \lambda_n\}$ are eigenvalues of ${\bf H}S{\bf H}$  in descending order, $i=1,\cdots, n$. In particular, when $\mu= {\bm 0}_m$, $\xi_i\stackrel{\text{iid} }\thicksim \chi^2_{m}$, $i=1,\cdots, n$.
\end{theorem}

This theorem shows that  the pseudo-$F$ test statistic  has the same asymptotic distribution as the random variable of Chi-squared-type mixtures. From the proofs presented in the Appendix, we can see that the asymptotic results do not depend on the dimension of response variables and are  valid for high-dimensional data. To derive the asymptotic distribution of the proposed test statistic $T_{sqrt}$, we need the following assumption.

\begin{assumption}\label{Assump1}
	$ \sum_{i=1}^{\infty} (\lambda^*_i)^{1/2} < \infty $ and  $\sum_{i=1}^{\infty} \big| (\lambda_i/n)^{1/2}-(\lambda^*_i)^{1/2} \big| \stackrel{P}\longrightarrow 0$.
\end{assumption}

Note that \citet{Gretton2009} also assumed that  $ \sum_{i=1}^{\infty} (\lambda^*_i)^{1/2} < \infty $ and  \citet{Zhang2012} shows that $(\lambda_i/n)^{1/2}\stackrel{P}\longrightarrow(\lambda^*_i)^{1/2}$ for $i=1,2,\cdots$. It implies that the assumption can be easily satisfied when
the number of eigenvalues of kernel $\psi(\cdot,\cdot)$ are finite, such as linear kernel. By the Theorem 1 of \citet{Gretton2009} and Assumption \ref{Assump1}, we have the following proposition.

\begin{proposition}\label{prop1}
	Assume that  $\{z_1, z_2,\cdots\}$ is an infinite sequence of {\it i.i.d.} standard Gaussian variables and Assumption \ref{Assump1} holds. Then $\sum_{i=1}^{\infty}\big((\lambda_i/n)^{1/2}-(\lambda^*_i)^{1/2}\big) z_i^2 \stackrel{p} \longrightarrow 0$ as  $n \rightarrow \infty$.
\end{proposition}

\begin{lemma}\label{lem3}
	For any symmetric semidefinite matrices $\bf D_1$ and $\bf D_2$ of  n dimension, the eigenvalues of ${\bf D}= {\bf D}_1 {\bf D}_2$  are nonnegative such that $\text{tr}({\bf D}^\top{\bf D}) \leq  \text{tr}({\bf D})^2 $.
\end{lemma}
When $\bf D$ is a symmetric semidefinite matrix, the result in \ref{lem3} naturally holds, since the sum of squares of the nonnegative eigenvalues should be smaller than the squares of the sum of them. 
Lemma \ref{lem3} is vital to show that a matrix formed by the product of  two symmetric semidefinite matrices still have the same property. Based on it and  proposition \ref{prop1}, the asymptotic null distribution of the proposed test statistic thus follows.

\begin{theorem}\label{thm2}
	When Assumptions $1$ holds, $T_{sqrt}$ has the same asymptotic distribution with $m^{-1}\sum_{i=1}^{n} \eta_i\xi_{i}$ when $n\rightarrow\infty$, where $\xi_i=\varpi_i+\varrho_i/(1+\mu^\top \Delta^{-1} \mu) $, $\varpi_i\stackrel{\text{iid} }\thicksim \chi^2_{m-1}$ and $\varrho_i\stackrel{\text{iid} }\thicksim \chi^2_1$
	are mutually independent, and $\eta_i=\lambda_i^{1/2}/\sum_{j=1}^{n} \lambda_j^{1/2}$, and $\{\lambda_1,\cdots, \lambda_n\}$ are eigenvalues of ${\bf H}S{\bf H}$  in descending order, $i=1,\cdots,n$. In particular, when $\mu= \bm 0_m$, $\xi_i\stackrel{\text{iid} }\thicksim \chi^2_{m}$, $i=1,\cdots,n$.
\end{theorem}

\subsection{Computation issue}
$T_{pseudo}$ and $T_{sqrt}$ have  asymptotical null distribution  of the same form with two mixtures of Chi-squared random variables, whose density functions involve multiple integrations. 
To ease computation and increase efficiency, we employ a parameter bootstrap procedure to
approximate them in finite-sample cases as the explicit formulas of the distributions of $T_{pseudo}$, $T_{sqrt}$ and the two mixtures are all intricate. 

\vspace{0.1cm}
\noindent{\bf A parameter bootstrap procedure } \\
{\bf\textit{Step 1.}} Set a large  $B$, for example $B=1000$.\\
{\bf\textit{Step 2.}} Randomly generate $n$ pairs of observations from the Chi-squared distributions with $m-1$ and $1$ degrees of freedom, denoted by $\{(\varpi_{i}, \varrho_i)|i=1,\cdots,n\}$. Let $T_1=m^{-1}\sum_{i=1}^{n} w_i [\varpi_i+\varrho_i/(1+\mu^\top \Delta^{-1} \mu)]$ and $T_2=m^{-1}\sum_{i=1}^{n} \eta_i [\varpi_i+\varrho_i/(1+ \mu^\top \Delta^{-1} \mu)]$.\\
{\bf\textit{Step 3.}} Repeat Step 2  $B$ times and denote the obtained   statistics $T_1$ and $T_2$ by $T_{11},\cdots,T_{1B}$ and $T_{21},\cdots,T_{2B}$.\\
{\bf\textit{Step 4.}} By the law of large numbers, the p-values of $T_{pseudo}$ and $T_{sqrt}$ can be empirically estimated by
$$p_{pseudo}=\frac{1}{B}\sum\limits_{i=1}^B I\{T_{1i}\geq t_{pseudo}\}$$and$$p_{sqrt}=\frac{1}{B}\sum\limits_{i=1}^B I\{T_{2i}\geq t_{sqrt}\},$$
where $t_{pseudo}$ and $t_{sqrt}$ is the observed value of $T_{pseudo}$ and $T_{sqrt}$ respectively and $I\{\cdot\}$ is an indicator function.

It is worth pointing out that  the proposed parameter bootstrap procedure is much faster than the original bootstrap procedure. For a sample of size $n$, the bootstrap procedure needs to generate $n$ individual observations,
perform $n(n-1)/2$  calculations to obtain  the similarity matrix and conduct two matrix multiplications, while the proposed procedure just needs to generate $2n$ random samples and calculate some  summations, which runs ver fast. For example, it takes 0.57 seconds to implement a parameter bootstrap for $B=1000$ and $n=500$ using Inter Core (TM) i9-9900 CPU, and takes 246.85 seconds to implement the original bootstrap procedure under the same setting.

In addition to the parameter bootstrap procedure, one can use the Box scaled $\chi^2$-approximation \citep{Box1954} or the generalized gamma distribution approximation \citep{Li2014} by matching several cumulants of them. Based on Theorems 1 and 2, the  cumulants of $T_{pseudo}$ and $T_{sqrt}$ can be estimated easily. In particular,  the $l$th cumulant of $T_{pseudo}$ and $T_{sqrt}$ are estimated as  $b_l= 2^{l-1}(l-1)! m_0\sum_{i=1}^{n} (w_i/m)^l$ and $h_l= 2^{l-1}(l-1)! m_0\sum_{i=1}^{n} (\eta_i/m)^l$ respectively, $l=1,2,3,4$, where $m_0=m-1+1/(1+\mu^\top \Delta^{-1} \mu) $, $w_i=\lambda_i/\sum_{j=1}^{n} \lambda_j$, $\eta_i=\lambda_i^{1/2}/\sum_{j=1}^{n} \lambda_j^{1/2}$,  and $\{\lambda_1,\cdots, \lambda_n\}$ are eigenvalues of ${\bf H}S{\bf H}$  in a descending order, $i=1,\cdots,n$. Let a random variable $X_g$ follow a generalized gamma distribution with the density function being
$$f_g(x)=\frac{v_gx^{v_g w_g-1}}{\sigma_g^{v_gw_g}\Gamma(w_g)}\exp\{-(\frac{x}{\sigma_g})^v\},x>0,$$
where $v_g$, $w_g$ and $\sigma_g$ are the parameters. We recommend to use the distribution  of  $X_g+\theta_g$ to approximate those of $T_{pseudo}$ and $T_{sqrt}$, where $\theta_g$ is a location parameter. The unknown parameters $v_g$, $w_g$, $\sigma_g$, and $\theta_g$ can be obtained by solving the equations
$$\left\{\begin{array}{ccl}
b_1~(\hbox{or}~h_1)&=&m_1+\theta_g\\
b_2~(\hbox{or}~h_2)&=&m_2-m_1^2\\
b_3~(\hbox{or}~h_3)&=&m_3-3m_2m_1+2m_1^3\\
b_4~(\hbox{or}~h_4)&=&m_4-4m_3m_1-3m_2^2+12m_2m_1^2-6m_1^4,
\end{array}\right.$$
where $m_l=E(X_g^l)=\sigma_g^l \Gamma(v_g^{-1}l+w_g)/\Gamma(w_g)$ is the $l$th moment of $X_g$, $l=1,2,3,4$, and $\Gamma(\cdot)$ is a Gamma function.

\subsection{Power analysis  }
To study the asymptotic power of the pseudo-$F$ test, we consider the multivariate linear model $Y=X \bm \beta + \varepsilon$, where $Y=(y_1,\ldots,y_n)^\top$ with $y_i=(y_{i1},\cdots, y_{iq})^\top$, ${\boldsymbol\beta}=(\beta_{ij})_{m\times q}$ is the effect of the predictor variables  $\mathbb X$ on the responses variables  $\mathbb Y$ and $\varepsilon=(\varepsilon_1,\cdots, \varepsilon_n)^\top$ are random errors with $E(\varepsilon_i)={\bf0}_{n\times q}$ and $\text{cov}(\varepsilon_i)=\Delta_{\varepsilon}$ for $i=1,\cdots,n$.  Then for a linear kernel $\psi$,   i.e., $s_{ij}=\psi(\bm y_i,\bm y_j)$, we have the following theorem.

\begin{theorem}\label{thm3}
	Assume that for $0<\iota<0.5$, $\max_{1\leq  i\leq  m,1\leq  j\leq  q} n^{0.5-\iota}|\beta_{ij}|=\tilde{c} \neq 0$ as $ n\rightarrow \infty$, then for the linear kernel $\psi(\cdot,\cdot)$ and $\tilde\tau_0>0$, we have
	$$\lim_{n\rightarrow\infty}  P\Big(\Big|\frac{T_{pseudo}}{n^{2\iota}}- \frac 1m \frac{\text{tr}(\tilde\Delta^{-1} \Delta \tilde{\boldsymbol{\beta}} \tilde{\boldsymbol{\beta}}^\top \Delta^\top ) }{\text{tr}(\Delta_{\varepsilon})}\Big|>\tilde\tau_0\Big)=0.$$
	It thus follows that $P\big(T_{pseudo}> F_1^{-1}(1-\alpha)\big)\stackrel{p} \longrightarrow 1$ as $n\rightarrow\infty$, where $F_1$ is the asymptotic distribution function of $T_{pseudo}$, $F_1^{-1}(1-\alpha)$ is the $(1-\alpha)$-quantile of $F_1$ and $\alpha$ is the nominal significance level.
\end{theorem}

Theorem \ref{thm3} shows that when the effect $\bm \beta$ is relatively small, even as small as $n^{\iota-0.5}$ with  tiny $\iota$, the power of the pseudo-$F$ test still converges to one in probability as $n\rightarrow\infty$. Actually, when $\bm \beta$ is fixed, the asymptotic power of the square-root $F$-type test statistic has the same performance.

\begin{theorem}\label{thm4}	
	For the linear kernel $\psi(\cdot,\cdot)$, $P\big(T_{sqrt}> F_2^{-1}(1-\alpha)\big)\stackrel{p} \longrightarrow 1$ as $n\rightarrow\infty$, where $F_2$ is the asymptotic distribution function of $T_{sqrt}$, $F_2^{-1}(1-\alpha)$ is the $(1-\alpha)$-quantile of $F_2$, and $\alpha$ is the nominal significance level.
\end{theorem}

\section{Simulation studies}\label{sim}

\subsection{Accuracy of p-value calculation}

We evaluate through simulation studies the p-value calculation accuracy of the asymptotic distributions of the pseudo-$F$ ($T_{pseudo}$) and square-root $F$-type  ($T_{sqrt}$) test statistics for various correlation matrices. The similarity matrix for the response variables is constructed based on the measure of inner product, that is,  $S=YY^\top$. The observations for the predictor variables $\mathbb X$ are generated from the $m$-dimensional normal distribution $N({\bf1}_m, {\it \boldsymbol\Theta}_x)$ with mean vector ${\bf1}_m=(1,\cdots, 1)^\top$ and covariance matrix ${\it \boldsymbol\Theta}_x=\big(\theta_{ij}^{(x)}\big)_{m\times m}$, where $\theta_{ij}^{(x)}=\rho_x^{|i-j|}$ with $\rho_x=0.5$. The observations for  $\mathbb Y$ are generated from $N_k({\bf0}_k, {\it \boldsymbol\Theta}_y)$. We consider the following two correlation models for ${\it \boldsymbol\Theta}_y=\big(\theta_{ij}^{(y)}\big)_{k\times k}$.

\begin{itemize}
	\item Model 1 (AR(1) correlation): Let $\theta_{ij}^{(y)}=\rho_y^{|i-j|}$ for $1\leq i,j\leq k$, where $\rho_y=0.3, 0.8$.
	
	\item Model 2 (Equal correlation): Let $\theta_{ij}^{(y)}=\rho_y$ for $1\leq i\neq j\leq k$ and $\theta_{ii}^{(y)}=1$ for $1\leq i\leq k$, where $\rho_y=0.3, 0.8$.
\end{itemize}
Thus there are a total of four correlation matrices for the outcome variables $\mathbb Y$. We set the sample size $n$ to be 500. Let $m=5$ and $k=10$. For each correlation setting, 10000 simulation replicates are performed to  evaluate the empirical size of the tests at significance levels ranging from 0 to 1. In each simulation, the p-values of the $T_{pseudo}$ and $T_{sqrt}$ are calculated based on the asymptotic distribution and Monte Carlo method with $B=2000$ replicates.

Figures \ref{tab:fig1} and \ref{tab:fig2} below display the empirical sizes of $T_{pseudo}$ and $T_{sqrt}$  based on the asymptotic distribution against  significance levels under various correlation matrix settings. The figures show that the empirical sizes of $T_{pseudo}$ and $T_{sqrt}$ are always very close to the corresponding significance levels, even in the case of small significance levels. It indicates good accuracy of the asymptotic distributions derived in Theorems \ref{thm1} and \ref{thm2} and such accuracy  is not sensitive to the correlation structure and magnitude.

\begin{figure}
	%\captionsetup{justification=normal}
	\begin{center}
		\includegraphics[scale=0.8]{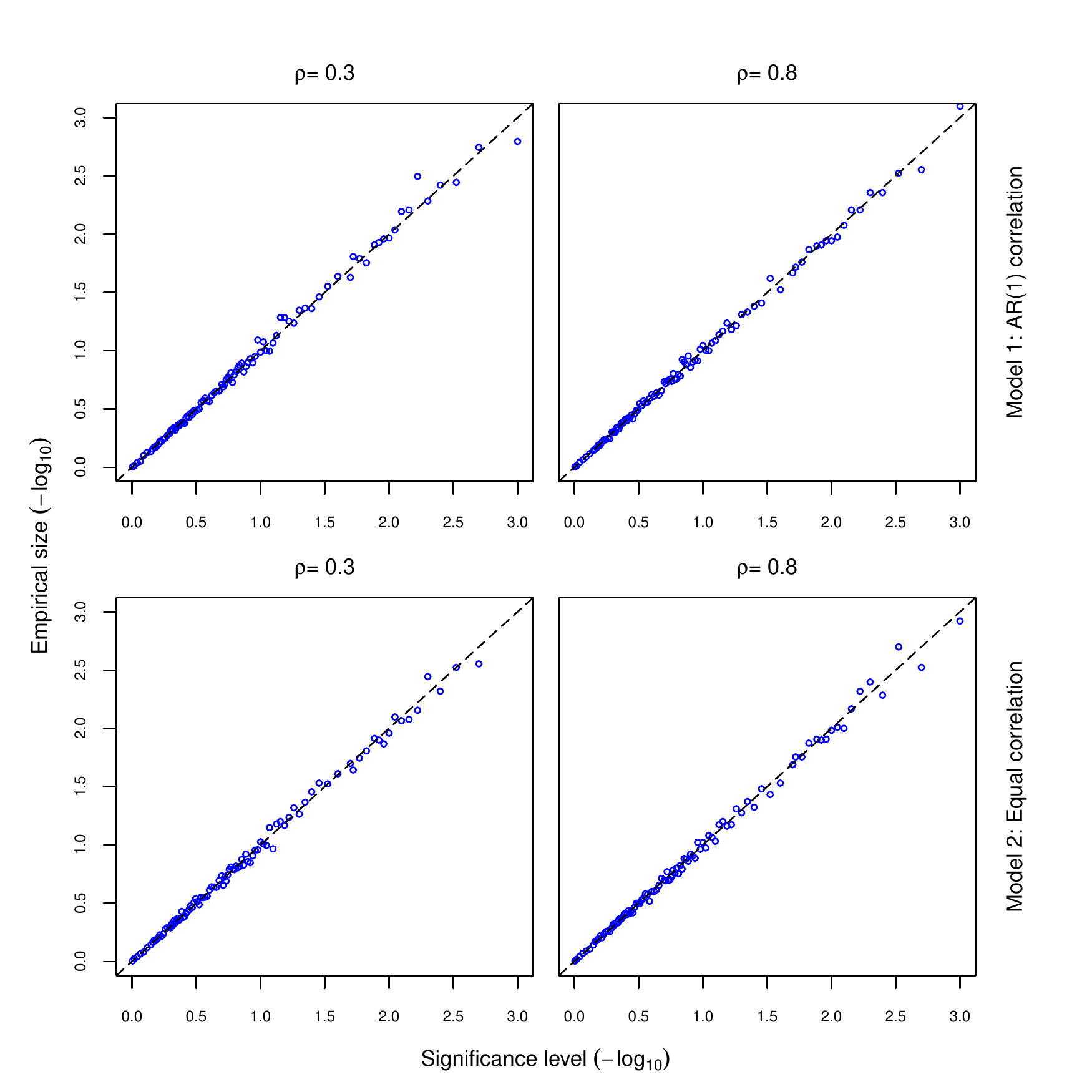}
	\end{center}\vspace{-0.5cm}
	\caption{Empirical sizes ($-\log_{10}$) of  the pseudo-$F$ test statistic ($T_{pseudo}$) based on the asymptotic distribution against the significance levels ($-\log_{10}$) for two correlation models.}
	\label{tab:fig1}
\end{figure}

\begin{figure}
	%\captionsetup{justification=normal}
	\begin{center}
		\includegraphics[scale=0.8]{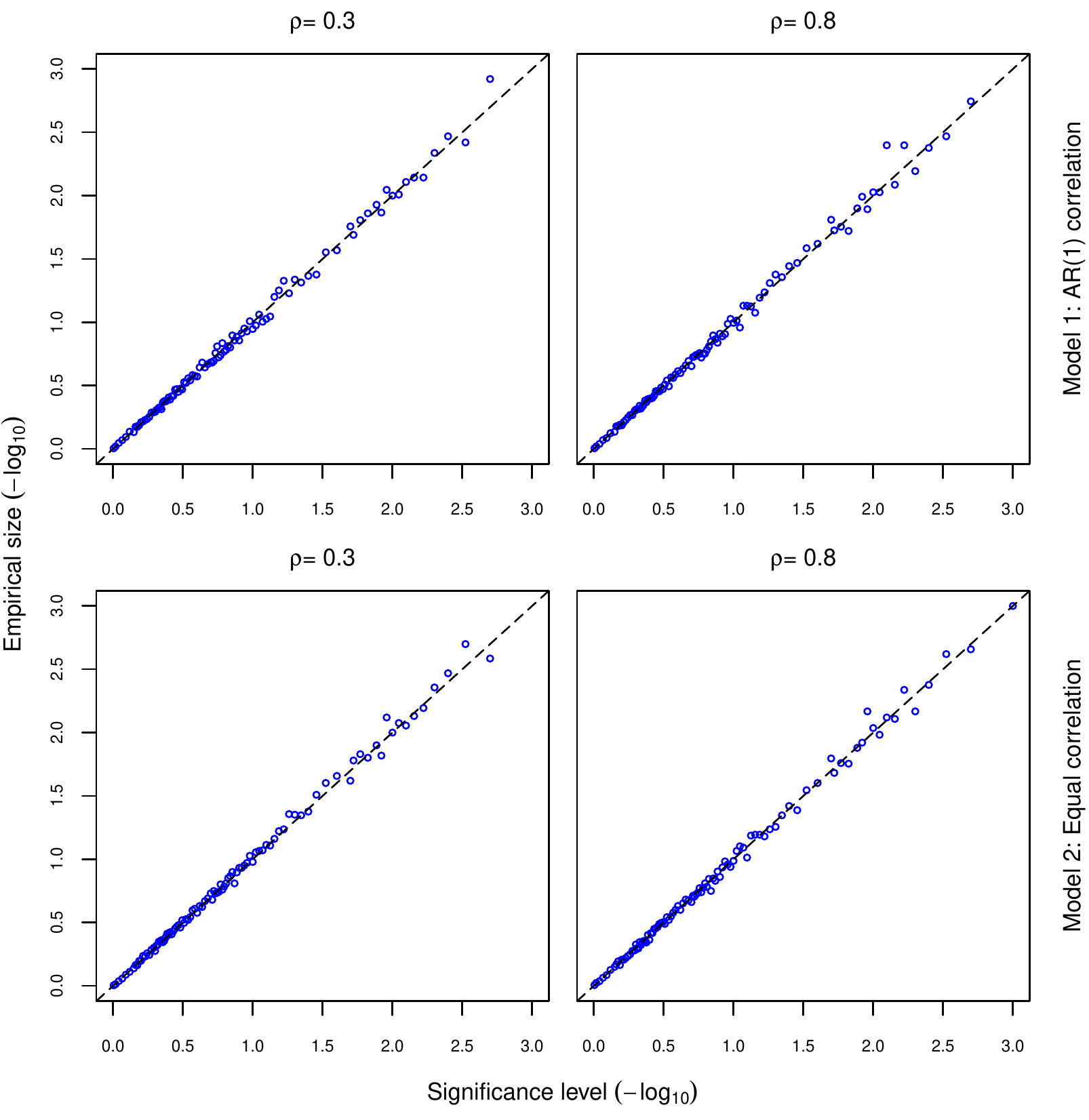}
	\end{center}\vspace{-0.5cm}
	\caption{Empirical sizes ($-\log_{10}$) of  the proposed  square-root $F$-type test statistic  ($T_{sqrt}$) based on the asymptotic distribution against the significance levels ($-\log_{10}$) for two correlation models.}
	\label{tab:fig2}
\end{figure}

\subsection{Power comparison}
Next we compare the type I error rates and powers of  $T_{pseudo}$ and $T_{sqrt}$. The multivariate response data are generated based on the linear model
$Y=X {\it \boldsymbol\beta} + \varepsilon$, where $\varepsilon\sim N_k({\bf0}_k, {\it \boldsymbol\Theta}_{\varepsilon})$.  The correlation matrix ${\it \boldsymbol\Theta}_{\varepsilon}=\big(\theta_{ij}^{(\varepsilon)}\big)_{k\times k}$ is set to have  the following two  correlation structures.

\begin{itemize}
	\item Model 1 (AR(1) correlation): Let $\theta_{ij}^{(\varepsilon)}=\rho_{\varepsilon}^{|i-j|}$ for $1\leq i,j\leq k$, where $\rho_{\varepsilon}=0.3, 0.8$.
	
	\item Model 2 (Equal correlation): Let $\theta_{ij}^{(\varepsilon)}=\rho_{\varepsilon}$ for $1\leq i \neq j\leq k$ and $\sigma_{ii}^{(\varepsilon)}=1$ for $1\leq i \leq k$, where $\rho_{\varepsilon}=0.3, 0.8$.
\end{itemize}
The observations $X$ for the predictor variables are generated from the multivariate normal distribution $N({\bf1}_m, {\it \boldsymbol\Theta}_x)$ with ${\it \boldsymbol\Theta}_x=\big(\theta_{ij}^{(x)}\big)_{m\times m}$ and $\theta_{ij}^{(x)}=0.5^{|i-j|}$. For the two tests, we define the similarity matrix as $S=YY^\top$. To investigate the performance of the two tests under different sparsity of the ``signals'', we choose the percentage $\tau$ of nonzero elements of ${\it \boldsymbol\beta}$ from $\{0\%, 20\%, 40\%, 60\%, 80\%, 100\%\}$. The null hypothesis corresponds to the case  $\tau =0\%$.  The signals (i.e., nonzero elements of ${\it \boldsymbol\beta}$) are set to be
$\{\log(k)/(25\tau km)\}^{1/2} + (1/k)\times N(0, 0.01)$, to make the powers comparable in different settings, where $N(0, 0.01)$ is the normal distribution with mean 0 and variance 0.01 and $\tau km$ is the total number of signals. We set $k=10, m=5$ and $n=500$. The empirical type I error rates and powers of the tests are calculated based on 1000 simulation replicates; in each simulation, 2000 Monte Carlo samples are drawn to calculate the p-values of the tests based on the asymptotic distributions. The nominal significance level is  0.05.

\begin{table}
	\caption{Type I error rates and powers of the pseudo-$F$ ($T_{pseudo}$) and square-root F-type   ($T_{sqrt}$) test statistics under two correlation models.}\label{tab:table1}
	\begin{center}
		%\addtolength{\tabcolsep}{-3pt}
		\begin{tabular}{cccccccc}
			\hline
			Correlation& Precentage  &&\multicolumn{2}{c}{$\rho=0.3$}&&\multicolumn{2}{c}{$\rho=0.8$}\\
			\cline{4-5}\cline{7-8}
			Model& of nonzeros&&$T_{pseudo}$ & $T_{sqrt}$ && $T_{pseudo}$ & $T_{sqrt}$ \\
			\hline
			1 &0\% &&0.057 &0.054&&0.051&0.048\\
			&20\% &&0.909 &0.886&&0.507&0.940\\
			&40\% &&0.825 &0.741&&0.395&0.663\\
			&60\% &&0.874 &0.794&&0.463&0.603\\
			&80\% &&0.884 &0.778&&0.478&0.597\\
			&100\%&&0.881 &0.803&&0.469&0.493\\
			\hline
			2 &0\%  &&0.051 &0.053&&0.051&0.049\\
			&20\%  &&0.787 &0.957&&0.306&1.000\\
			&40\%  &&0.669 &0.810&&0.249&0.980\\
			&60\%  &&0.706 &0.778&&0.280&0.930\\
			&80\%  &&0.683 &0.652&&0.305&0.718\\
			&100\% &&0.710 &0.571&&0.305&0.390\\
			\hline
		\end{tabular}
	\end{center}
\end{table}

Table \ref{tab:table1} presents the empirical type I error rates and powers of the two tests for various ${\boldsymbol\Sigma}_x$ and $\tau$. It can be seen from the table that both tests can control the type I error rates adequately and the proposed test $T_{sqrt}$  is generally more powerful than $T_{pseudo}$, especially when the correlation coefficients of the response variables are large. For example, under the AR(1) correlation model with $\rho=0.8$, the powers of $T_{pseudo}$ and $T_{sqrt}$ for $\tau=20\%$ are 0.507 and 0.940.  The superiority of $T_{sqrt}$ diminishes as the sparsity level of signals increases, but can still outperform $T_{pseudo}$ when the correlation coefficient is large. This implies that the proposed test tends to gain more power than $T_{pseudo}$ when the signals are sparse. $T_{pseudo}$ is slightly more powerful than $T_{sqrt}$ when the response variables are weakly dependent with decaying correlation (e.g., AR(1) correlation model with $\rho=0.3$). For example, under the AR(1) correlation model with $\rho=0.3$, the powers of $T_{pseudo}$ and $T_{sqrt}$ for $\tau=60\%$ are 0.874 and 0.794. In summary, the proposed test  has more robust power than the pseudo-$F$ statistic with respect to various sparsity levels of signals and correlation magnitudes; it performs consistently well across all settings.

\section{Applications}

We exemplify the tests using  the gene expression data from a  prostate cancer study \citep{Singh2002}, whose aim is to investigate whether gene expression differences underlie common clinical and pathological features of prostate cancer. This can be achieved by comparing gene expression differences between tumor and normal prostate samples. The data of expression profiles of approximately 12,600 genes from 52 tumor and 50 normal prostate specimens were collected. We confine our analysis to the genes of 10 pathways: Alanine, aspartate and glutamate metabolism (map00250); Pathogenic Escherichia coli infection (map05130); Viral myocarditis (map05416); PPAR signaling pathway (map03320); Rheumatoid arthritis (map05323); Tight junction (map04530); Regulation of actin cytoskeleton (map04810); Hypertrophic cardiomyopathy (map05410); Cardiac muscle contraction (map04260); TGF-beta signaling pathway (map04350). Detailed descriptions of these pathways are available in \cite{Wu1998}, \cite{Zihni2016}, \cite{Pinaud2018}, and among others.

Here for illustration, we consider using the gene expression patterns as the response variables and the status of tumor as the predictor variable. Likewise, the measure of inner product is used to construct the similarity matrix, that is, $S=YY^\top$. We are interested in whether the expression patterns of the genes in a pathway are associated with prostate tumor.  Table \ref{tab:table2} presents the p values of the pseudo $F$  and square-root $F$-type test statistics for testing the association between gene expression patterns and prostate tumor. The p value results are calculated  based on 10000 Monte Carlo replicates. This table shows that,  the p values   of the pseudo $F$ test are always larger than those of the square-root $F$-type test for the 10 pathways, indicating that the proposed test is more powerful than the pseudo F test. Moreover,  for some pathways such as map05130 and map04810, the pseudo $F$ test fails to detect any expression pattern difference between the normal and tumor samples under the significance level of 0.05, while the proposed test can  detect a difference.

\newpage
\renewcommand{\baselinestretch}{1.2}
\begin{table}
	\caption{P values of the pseudo F ($T_{pseudo}$) and square-root F-type ($T_{sqrt}$) test statistics for the association between gene expression patterns and prostate tumor.}\label{tab:table2}
	\begin{center}
		\begin{tabular}{cccc}
			\hline
			Pathway&Numbers of genes&$T_{pseudo}$ & $T_{sqrt}$\\
			\hline
			map00250& 33&0.0396 &0.0130\\
			map05130& 80&0.1797 &0.0134\\
			map05416& 98&0.0738 &0.0155\\
			map03320& 66&0.0055 &0.0025\\
			map05323&118&0.0097 &0.0010\\
			map04530&160&0.0564 &0.0072\\
			map04810&279&0.2297 &0.0121\\
			map05410&118&0.2385 &0.0185\\
			map04260& 72&0.0007 &0.0002\\
			map04350&128&0.0478 &0.0034\\
			\hline
		\end{tabular}
	\end{center}
\end{table}
\section{Conclusion}\label{disc}
Distance-based regression model is very effective to detect the relationship between  high-dimensional response variables and predictor variables of interest. It has a wide applications in many research fields related to statistics. One drawback is the intensive computation for the calculation of statistical significance, due to the lack of asymptotic null distribution. In this work, we establish the asymptotic distribution for the pseudo-$F$ test statistic based on the distance-based regression model and propose a  new test
which is more powerful than the original pseudo-$F$ test when the correlation coefficient of the outcomes for measuring the similarity is relatively large. The proposed theory is anticipated to further broaden the application of distance-based regression model.

The asymptotic property of the pseudo-$F$ test only requires that the kernel for the similarity matrix is a positive definite and does not impose any restrictions on the dimension of outcomes. So it is free to handle high-dimensional data and is expected to have a wide application in high-dimensional association studies.

\section*{Acknowledgement}
 We would like to thank the editor and two anonymous reviewers for their insighful comments.

\section*{Appendix}

%%%%%%%%  proof of Lemma \ref{lem1}  %%%%%%%%%
\subsection*{Proof of Lemma \ref{lem1}.}
Denote ${\it \Omega}= X\tilde\Delta^{-1}X^\top$. Note that ${\it\Omega}$ can be regarded as a random variable constructed based on a weighted inner product, which is a positive definite kernel. Then the numerator of $T_{pseudo}$ can be written as
\begin{equation*}
\text{tr}\big(H_X\tilde{S}\big)=\frac{1}{n}\text{tr}\big({\it\Omega} \tilde{S}\big)+ \text{tr}\big(({H_X-\frac{1}{n}{\it\Omega}} ) \tilde{S}\big).
\end{equation*}

We first show that  $\frac{1}{n}\text{tr}\big({{\it\Omega} {\bf H }S{\bf H } }\big)$ has the same asymptotical distribution with  $\frac{1}{n}\sum_{i=1}^{n} \lambda_i \xi_{i}$  under $H_0$. To this end, we denote  the eigenvalues of $X\tilde\Delta^{-1}X^\top$ by $\tilde\lambda_1 \geq   \cdots \geq  \tilde\lambda_m\geq  0$ and $\tilde\lambda_{m+1}=\cdots=\tilde\lambda_n =0$.
Let $\tilde\lambda^*_1\geq   \cdots \geq  \tilde\lambda ^*_n \geq  0$ be the solutions to  the  equation
\begin{equation}\label{App-Eq1}
\int \tilde\psi (x_1,x_2) \tilde u_i (x_1)p(x_1)dx_1=\tilde\lambda^*_i \tilde u_i(x_2)
\end{equation}
with $\tilde\psi(x_1,x_2)= x_1^\top\tilde\Delta^{-1} x_2-E_{x_1}( x_1^\top\tilde\Delta^{-1} x_2)-E_{x_2}( x_1^\top\tilde\Delta^{-1} x_2)+E_{x_1,x_2}( x_1^\top\tilde\Delta^{-1} x_2)=(x_1-\mu)^\top\tilde\Delta^{-1} (x_2- \mu),$
where $\tilde u_i(x_1)$ is the eigenfunction of the kernel $\tilde\psi(x_1,x_2)$ corresponding to $\tilde\lambda^*_i $. Note that eigenfunctions are orthonormal with respect to the probability measure $p(x_1)$, i.e.,
\begin{equation}\label{App-Eq2}
\int \tilde u_i (x_1) \tilde u_j (x_1)p(x_1)dx_1=\mathbb{I}_{i=j},
\end{equation}
where $\mathbb{I}_{\bf E} $ denoting the indicator function of ${\bf E}$. \citet{Zhang2012} shows that $\frac 1n \tilde\lambda_{i} \stackrel{p}\longrightarrow \tilde\lambda^*_{i}$,for $i=1,\cdots,n$.
Then  $\tilde\lambda^*_{m+1}= \cdots = \tilde\lambda^*_n=0$.

When $\mu=\bm 0$, it then follows that
$E_{x_1}(e_i^\top \Delta^{-1/2}{x_1 x_1}^\top \tilde\Delta^{-1} x_2 )= e_i^\top \Delta^{-1/2}x_2 $ and \\
$ E_{x_1}(e_i^\top \Delta^{-1/2}{x_1 x_1}^\top \Delta^{-1/2} e_j )=\mathbb{I}_{i=j} $,
where $\Delta^{-1/2}\Delta^{-1/2}=\Delta$. This implies that  $ \tilde\lambda^*_i =1$ and  $\tilde u_i(x_1)= e_i^\top \Delta^{-1/2} x_1$ for $i=1,\cdots,m$. When $\mu \neq \bm 0$, for $i,j=1,\cdots,m-1,$ it can be proved that $ \tilde\lambda^*_i =1$ and $ \tilde\lambda^*_m =1/(1+\mu^\top \Delta^{-1} \mu)$ corresponding to $\tilde u_i(x_1)= v_i^\top (x_1-\mu)$ and $\tilde u_m(x_1)= v_m^\top (x_1-\mu)$  satisfy  the equations \eqref{App-Eq1} and \eqref{App-Eq2}, where $v_m=  \Delta^{-1} \mu/ \mu^\top \Delta^{-1} \mu $ and $v_1,\cdots, v_{m-1}$ are the solutions to the equations $ v_i^\top \mu=0$ and $v_i^\top \Delta v_j=\mathbb{I}_{i=j}$.

We now show that  $ \tilde\lambda^*_1=\cdots=\tilde\lambda^*_{m-1} =1$, $ \tilde\lambda^*_m =1/(1+\mu^\top \Delta^{-1} \mu)$ and
$ \tilde\lambda^*_{m+1}=\cdots=\tilde\lambda^*_{n} =0$. By Theorem 3 in \citet{Zhang2012}, under $H_0$, $\frac 1n \text{tr}\big({ {\it\Omega} {\bf H }S{\bf H } }\big)$ has the same asymptotic distribution as $\sum_{i,j=1}^{n} \lambda^*_i \tilde\lambda^*_j z^2_{ij}=\sum_{i=1}^{n}  \lambda^*_i  \xi_i$, where $z_{ij}$ are {\it iid} standard Gaussian variables and $\xi_i=\sum_{j=1}^{m}\tilde\lambda^*_j z^2_{ij}=\varpi_i+\varrho_i/(1+\mu^\top \Delta^{-1} \mu) $ . In addition, by Theorem 1 in \citet{Gretton2009}, we have $\sum_{i=1}^{\infty} (\frac 1n \lambda_i - \lambda^*_i ) \xi_i \stackrel{p} {\longrightarrow} 0 $.  That is, $\sum_{i=1}^{n} \frac 1n \lambda_i \xi_i $  has the same asymptotic distribution as $\sum_{i=1}^{n}  \lambda^*_i  \xi_i $. It thus follows that $\frac{1}{n} \text{tr}\big({ {\it\Omega} {\bf H }S{\bf H } }\big)$ has the same asymptotic distribution as $\frac{1}{n}\sum_{i=1}^{n} \lambda_i \xi_{i}$.

Next we show that $\text{tr}\big(({H_X-\frac 1n {\it\Omega}} ) \tilde{S}\big)$ converges to 0 in probability.
Write
$$\text{tr}\big(({H_X-\frac 1n {\it\Omega}} ) \tilde{S}\big)=\text{tr}(AB)=\sum_{i=1}^{m}\sum_{j=1}^{m} a_{ij} b_{ij},$$
where  $A=(a_{ij})_{m \times m}= (\frac{1}{n}X^\top X)^{-1}-\tilde\Delta^{-1} $ and  $B=(b_{ij})_{m \times m}= \frac{1}{n}X^\top\tilde{S}X$.  By the law of large numbers, it can be obtained that the $(i,j) $th entry of $\frac{1}{n}X^\top X=\frac{1}{n} \sum_{i=1}^{n}x_i x_i^\top$ converges in probability to $\tilde\delta_{ij}$, that is
$a_{ij} \stackrel{p}\longrightarrow 0 $, $ i,j=1,\cdots,m$. Note that $(e_i-e_j)^\top S (e_i-e_j) \geq  0$ and $(e_i+e_j)^\top S (e_i+e_j) \geq  0$ for a positive definite kernel matrix $S$. Then we have $2|s_{ij}|\leq  s_{ii}+s_{jj}$ and $|E(s_{ij})|\leq  E(s_{ii})$, $i,j=1,\cdots,m$. It follows that
\begin{align*}
&E(s_{ij}^2) \leq  \frac{1}{4} E((s_{ii}+s_{jj})^2)  = \frac{1}{2} E(s_{ii}^2) +\frac{1}{2} E(s_{ii}s_{jj}) \leq   E(s_{ii}^2),\\
&|E(s_{ij}s_{lk})| \leq  \frac{1}{2} E(s_{ij}^2) +\frac{1}{2} E(s_{lm}^2) \leq   E(s_{ii}^2)=c_1,
\end{align*}
for  $i,j,l,k=1,\cdots,m$. Similarly, we can obtain that $|E(x_{j_1i_1} x_{j_2i_1}  x_{j_3i_2} x_{j_4i_2})|  \leq  \frac 14 E\big(x_{j_1i_1}^4+x_{j_2i_1}^4+x_{j_3i_2}^4 +x_{j_4i_2}^4\big) \leq   c_0$ for $i_1,i_2=1,\cdots,m$ and $j_1,j_2,j_3,j_4=1,\cdots,n$.

Define $ \bar x=\frac{1}{n} \sum_{j=1}^{n} x_j$ and $\bar x=(\bar{x}_{1}, \cdots,\bar{x}_{m})$. Then we can write
\begin{align*}
B=\frac{1}{n} \sum_{i=1}^{n}\sum_{j=1}^{n}(x_i-\bar x)s_{ij}(x_j-\bar x)^\top.
\end{align*}

Owing to $x_i-\bar x=(x_i-\mu)-(\bar x-\mu)$, we can assume that $E(\mathbb X )=\bm 0_m$ to estimate the expectation and variance of $b_{ij}$ under $H_0$,
\begin{equation*}
\begin{aligned}
|E(B)|=&\big|\frac{1}{n} \sum_{i=1}^{n}\sum_{j=1}^{n} E(s_{ij})E\big((x_i-\bar x)(x_j-\bar x)^\top \big)\big|\\
\leq&  E(s_{ii}) \frac{1}{n} \sum_{i=1}^{n} \big|E\big((x_i-\bar x)(x_i-\bar x)^\top \big)\big|+E(s_{ii}) \frac{1}{n} \sum_{i\neq  j=1}^{n}\big|E\big((x_i-\bar x)(x_j-\bar x)^\top \big)\big|\\
\leq&   E(s_{ii}) |\Delta|+ E(s_{ii}) |\Delta| =2E(s_{ii}) |\Delta|,
\end{aligned}
\end{equation*}
implying that $E(b_{i_1i_2})$ is finite, $i_1,i_2=1,\cdots,n$. Through some algebraic manipulations, it can be obtained that

\begin{align*}
E(b_{i_1 i_2}^2)=&\dfrac{1}{n^2} \sum_{j_1,j_2,j_3,j_4=1}^{n}
c_s
E\big(({ x_{j_1i_1}}-\bar{ x}_{i_1})({ x_{j_2i_1}}-\bar{ x}_{i_1}) ({ x_{j_3i_2}}-\bar{ x}_{i_2}))({ x_{j_4i_2}}-\bar{ x}_{i_2}) \big)\\
=&\dfrac{1}{n^2} \sum_{j_1,j_2,j_3,j_4=1}^{n}  c_s E\big
({ x_{j_1i_1}}{ x_{j_2i_1}} { x_{j_3i_2}}{ x_{j_4i_2}} \big)
-\dfrac{4}{n^2} \sum_{j_1,j_2,j_3,j_4=1}^{n}  c_s E\big
({ x_{j_1i_1}}{ x_{j_2i_1}} { x_{j_3i_2}}{\bar{ x}_{i_2}} \big)\\
&+\dfrac{2}{n^2} \sum_{j_1,j_2,j_3,j_4=1}^{n}  c_s E\big
({ x_{j_1i_1}}{ x_{j_2i_1}} {\bar{ x}_{i_2}}{\bar{ x}_{i_2}} \big)
+\dfrac{4}{n^2} \sum_{j_1,j_2,j_3,j_4=1}^{n}  c_s E\big
({ x_{j_1i_1}}{\bar{ x}_{i_1}}{ x_{j_3i_2}} {\bar{ x}_{i_2}} \big)\\
&-\dfrac{4}{n^2} \sum_{j_1,j_2,j_3,j_4=1}^{n}  c_s E\big
({ x_{j_1i_1}}{\bar{ x}_{i_1}}{\bar{ x}_{i_2}}{\bar{ x}_{i_2}} \big)
+\dfrac{1}{n^2} \sum_{j_1,j_2,j_3,j_4=1}^{n}  c_s E\big
({\bar{ x}_{i_1}}{\bar{ x}_{i_1}}{\bar{ x}_{i_2}}{\bar{ x}_{i_2}} \big)\\
\triangleq&\kappa_1+\kappa_2+\kappa_3+\kappa_4+\kappa_5+\kappa_6,
\end{align*}
where $c_s=E(s_{j_1j_2}s_{j_3j_4})$. Then using the assumptions $E(\mathbb X_1^4)\leq  c_0$ and $E(s_{11}^2)\leq  c_1$, we have
\begin{align*}
|\kappa_1| \leq & \dfrac{c_1}{n^2}  \sum_{j_1,j_2,j_3,j_4=1}^{n}  \big| E\big
( x_{j_1i_1} x_{j_2i_1}  x_{j_3i_2}x_{j_4i_2} \big) \big|\\
\leq & \dfrac{c_1 }{n} \max_{i_1,i_2} E({ x^2_{1i_1}}{ x^2_{1i_2}})+ c_1\max_{i_1,i_2} E({ x^2_{1i_1}}{ x^2_{2i_2}})
+2 c_1 \max_{i_1,i_2} |E({ x_{1i_1}}{ x_{2i_1}}{ x_{1i_2}}{ x_{2i_2}})|\\
\leq & (\dfrac{1}{n}c_0c_1+c_0c_1 +2c_0c_1  ) \leq  4c_0c_1.
\end{align*}

Similarly,
\begin{align*}
|\kappa_2| \leq & \dfrac{4}{n^2} \sum_{j_1,j_2,j_3,j_4=1}^{n} |c_s|
|E({ x_{j_1i_1}}{ x_{j_2i_1}} { x_{j_3i_2}}{\bar{ x}_{i_2}})|\\
\leq & \dfrac{4c_1}{n^2} \sum_{j_1,j_2,j_3,j_4=1}^{n}
|E({ x_{j_1i_1}}{ x_{j_2i_1}} { x_{j_3i_2}}{ x_{j_4i_2}} )|
\leq  16c_0c_1,\\
|\kappa_3| \leq & \dfrac{2}{n^2} \sum_{j_1,j_2,j_3,j_4=1}^{n} |c_s|
|E({ x_{j_1i_1}}{ x_{j_2i_1}} {\bar{ x}_{i_2}}{\bar{ x}_{i_2}})|\\
\leq &\dfrac{2c_1}{n^2} \sum_{j_1,j_2,j_3,j_4=1}^{n}
|E({ x_{j_1i_1}}{ x_{j_2i_1}} { x_{j_3i_2}}{ x_{j_4i_2}} )|\leq 8c_0c_1,\\
|\kappa_4| \leq & \dfrac{4}{n^2} \sum_{j_1,j_2,j_3,j_4=1}^{n} |c_s|
|E({ x_{j_1i_1}}{\bar{ x}_{i_1}}{ x_{j_3i_2}} {\bar{ x}_{i_2}})|\\
\leq  &\dfrac{4c_1}{n^2} \sum_{j_1,j_2,j_3,j_4=1}^{n}
|E({ x_{j_1i_1}}{ x_{j_2i_1}} { x_{j_3i_2}}{ x_{j_4i_2}} )|\leq 16c_0c_1,\\
|\kappa_5| \leq & \dfrac{4}{n^2} \sum_{j_1,j_2,j_3,j_4=1}^{n} |c_s|
|E({ x_{j_1i_1}}{\bar{ x}_{i_1}} {\bar{ x}_{i_2}} {\bar{ x}_{i_2}})|\\
\leq & \dfrac{4c_1}{n^2} \sum_{j_1,j_2,j_3,j_4=1}^{n}
|E({ x_{j_1i_1}}{ x_{j_2i_1}} { x_{j_3i_2}}{ x_{j_4i_2}} )|
\leq  16c_0c_1,\\
|\kappa_6| \leq & \dfrac{1}{n^2} \sum_{j_1,j_2,j_3,j_4=1}^{n} |c_s|
|E({\bar{ x}_{i_1}}{\bar{ x}_{i_1}} {\bar{ x}_{i_2}} {\bar{ x}_{i_2}})|\\
\leq & \dfrac{c_1}{n^2} \sum_{j_1,j_2,j_3,j_4=1}^{n}
|E({ x_{j_1i_1}}{ x_{j_2i_1}} { x_{j_3i_2}}{ x_{j_4i_2}} )|\leq   4c_0c_1.
\end{align*}
It thus follows that $E(b_{i_1 i_2}^2)\leq  4c_0c_2(1+4+6+4+1)=64c_0c_1$ for $i_1,i_2=1,\cdots,m$ and $j_1,j_2,j_3,j_4=1,\cdots,n$. Then $E(b_{i_1i_2})$ and $\text{var}(b_{i_1i_2})$ are finite. Further by the Chebyshev's inequality, $b_{i_1i_2}$ is bounded in probability. Therefore, $\text{tr}\big(({H_X-\frac 1n {\it\Omega}} ) \tilde{S}\big)$ converges to 0 in probability. This completes the proof.

\vspace{2cm}

%%%%%%%%  proof of Lemma \ref{lem1} %%%%%%%%%
\subsection*{Proof of Lemma \ref{lem2}.}
Let $\mathbb H$ be a reproducing kernel Hilbert space on $\mathbb X$, with a continuous feature mapping $\phi (x)$. With the positive definite kernel function $s_{ij}=\psi(y_i,y_j)=\langle\phi ({ y_i}),\phi ({y_j})\rangle_{\mathbb H} $, it can be obtained that
\begin{align*}
\frac {1}{n^2 }{\bf 1}^\top_n S {\bf 1}_n=\langle\frac 1n \sum_{i=1}^{n}\phi ({ y_i}),\frac 1n \sum_{j=1}^{n}\phi ({ y_j})\rangle_{\mathbb H}
\stackrel{p} \longrightarrow E(\langle\phi ({ y_i}),\phi ({y_j})\rangle_{\mathbb H} )=E(s_{ij}).
\end{align*}

By the law of large numbers, we have
$$  \frac 1n\sum_{i=1}^{n}\lambda_i=\frac 1n \text{tr}(\tilde{S})=\frac 1n \text{tr}(S)-\frac {1}{n^2 }{\bf 1}^\top_n S {\bf 1}_n
\stackrel{p} \longrightarrow E(s_{ii})-E(s_{ij}).$$
By Lemma \ref{lem1}, $\text{tr}\big(H_X\tilde{S}\big)/\sum_{i=1}^{n} \frac {\lambda_i}{ n} $ has the same asymptotic distribution as $\sum_{i=1}^{n} w_i \xi_{i}$, where $w_i=\lambda_i/\sum_{j=1}^{n} \lambda_j$. Since $\sum_{i=1}^{n} w_i=1$ and $0\leq  w_i \leq  1$, $i=1,\cdots,n$, it then follows that
$$E\big(\sum_{i=1}^{n} w_i \xi_{i}\big)=m_0\sum_{i=1}^{n} w_i=m_0~~ \text{and} ~~
\text{var}\big(\sum_{i=1}^{n} w_i \xi_{i}\big)=2m_0\sum_{i=1}^{n} w_i^2 \leq  2m_0\sum_{i=1}^{n} w_i =2m_0, $$
where $m_0=1/(1+\mu^\top \Delta^{-1} \mu)+m-1$. Hence, by Chebyshev's inequality, for any $\tau_0>0$,
\begin{align*}
P\Big(\Big|\frac 1n \sum_{i=1}^{n} w_i \xi_{i}\Big|\geq  \tau_0+\frac{ m_0}{n}  \Big)& \leq
P\Big(\Big|\frac 1n \sum_{i=1}^{n} w_i \xi_{i}-\frac {m_0}{n} \Big|\geq  \tau_0 \Big)\\ &\leq  \frac{2m_0}{n^2\tau_0^2}.
\end{align*}

Then for any $\tau_1=2\tau_0>0$,
$$\lim\limits_{n \rightarrow \infty} P\Big(\Big|\frac 1n \sum_{i=1}^{n} w_i \xi_{i}\Big|\geq  \tau_1 \Big)=0.$$
It follows that $\text{tr}\big(H_X\tilde{S}\big)/\sum_{i=1}^{n} \lambda_i$ converges in distribution to zero, that is,  $\text{tr}\big(H_X\tilde{S}\big)/\sum_{i=1}^{n} \lambda_i$ converges in probability to zero. Then we conclude that
\begin{align*}
\dfrac{1}{n-m}\text{tr}\big(({\bf I}_n - H_X)\tilde{S}\big)&=\frac{\text{tr}\big(({\bf I}_n - H_X)\tilde{S}\big) }{\sum_{i=1}^{n}\lambda_i } \dfrac{1}{n-m}\sum_{i=1}^{n}\lambda_i \stackrel{p} \longrightarrow E(s_{ii})-E(s_{ij}).
\end{align*}

%%%%%%%%  proof of Theorem \ref{thm1}  %%%%%%%%%
\subsection*{Proof of Theorem \ref{thm1}.} Theorem \ref{thm1} is a direct consequence of Lemma \ref{lem1} and Lemma \ref{lem2} .

%%%%%%%%  proof of Proposition \ref{prop1}  %%%%%%%%%
\subsection*{Proof of Proposition \ref{prop1}.}
When Assumption \ref{Assump1} holds, by Theorem 1  of \citet{Gretton2009} and the Chebyshev's inequality, we have  $ \sum_{i=1}^{\infty}(\frac{\lambda_i}{n})^{1/2} z_i^4$ and $ \sum_{i=1}^{\infty}(\lambda_i^*)^{1/2} z_i^4$ are bounded in probability, as $n \rightarrow \infty$. Combining this and the Cauchy-Schwarz inequality leads to
\begin{align*}
\Big| \sum_{i=1}^{\infty} \big((\frac{\lambda_i}{n})^{1/2}-(\lambda^*_i)^{1/2}\big) z_i^2 \Big|
\leq & \Big\{  \sum_{i=1}^{\infty}(\frac{\lambda_i}{n})^{1/2} z_i^4 \Big\}^{1/2} \Big\{   \sum_{i=1}^{\infty} \Big| (\frac{\lambda_i}{n})^{1/4}-(\lambda^*_i)^{1/4} \Big|^2 \Big\}^{1/2}\\
&+\Big\{  \sum_{i=1}^{\infty}(\frac{\lambda_i}{n})^{1/2} z_i^4 \Big\}^{1/2}
\Big\{   \sum_{i=1}^{\infty} \Big| (\frac{\lambda_i}{n})^{1/4}-(\lambda^*_i)^{1/4} \Big|^2 \Big\}^{1/2}\\
\leq & \Big\{  \sum_{i=1}^{\infty}(\frac{\lambda_i}{n})^{1/2} z_i^4 \Big\}^{1/2} \Big\{   \sum_{i=1}^{\infty} \Big| (\frac{\lambda_i}{n})^{1/2}-(\lambda^*_i)^{1/2} \Big| \Big\}^{1/2}\\
&+\Big\{  \sum_{i=1}^{\infty}(\lambda_i^*)^{1/2} z_i^4 \Big\}^{1/2}
\Big\{   \sum_{i=1}^{\infty} \Big| {(\frac{\lambda_i}{n})}^{1/2}-{\lambda^*_i}^{1/2} \Big| \Big\}^{1/2}
\stackrel{P}\longrightarrow 0.
\end{align*}

%%%%%%%%  proof of Lemma \ref{lem3} %%%%%%%%%

\subsection*{Proof of Lemma \ref{lem3}.}
For any symmetric semidefinite matrices ${\bf D}_1$, there exits a symmetric matrices $\bf P$, so that  ${\bf D}_1={\bf P} {\bf P}^\top$.
For any n dimension vector $\bf x$, we have ${\bf x} ^\top {\bf P}^\top {\bf D}_2 {\bf P} {\bf x} \geq  0$. It follows that ${\bf P}^\top {\bf D}_2 {\bf P}$ is a symmetric semidefinite matrice. Thus ${\bf D}={\bf D}_1 {\bf D}_2= {\bf PP}^\top {\bf D}_2$ has the same eigenvalues as ${\bf P}^\top {\bf D}_2 {\bf P}$ whose eigenvalues are nonnegative. Consequently, $\text{tr}({\bf D}^\top{\bf D}) =\sum_{i=1}^{n} \lambda_{i,D}^2\leq  \big(\sum_{i=1}^{n} \lambda_{i,D}\big)^2 =\text{tr}({\bf D})^2 $, where $\lambda_{i,D}$ is eigenvalue of $\bf D$, $i=1,\cdots,n$.

%%%%%%%%  proof of Theorem \ref{thm2} %%%%%%%%%
\subsection*{Proof of Theorem \ref{thm2}.}
For ${\it\Omega}= X\tilde\Delta^{-1}X^\top$, it can be obtained that
\begin{equation*}
n^{1/2}\text{tr}\big(H_X\tilde{S}^{1/2}\big)=\text{tr}\big({\it\Omega}(\tilde{S}/n)^{1/2}\big)+ n^{1/2}\text{tr}\big(({H_X-\frac{1}{n}{\it\Omega}} ) \tilde{S}^{1/2}\big).
\end{equation*}

We first show that  $\text{tr}\big({{\it\Omega} (\tilde{S}/n)^{1/2} }\big)$ has the same asymptotical distribution with  $\sum_{i=1}^{n} (\lambda_i/n)^{1/2} \xi_{i}$  under $H_0$.  By extending the proof of Theorem 3 in \citet{Zhang2012}, under $H_0$, $\text{tr}\big({{\it\Omega} (\tilde{S}/n)^{1/2} }\big)$ has the same asymptotic distribution as $\sum_{i,j=1}^{n} (\lambda^*_i)^{1/2} \tilde\lambda^*_j z^2_{ij}=\sum_{i=1}^{n}  (\lambda^*_i)^{1/2}  \xi_i$. In addition, by proposition \ref{prop1}, we have $\sum_{i=1}^{n} \big((\lambda_i/n)^{1/2} -(\lambda^*_i )^{1/2} \big)\xi_i \stackrel{p}\longrightarrow 0 $. It thus follows that  $\text{tr}\big({{\it\Omega} (\tilde{S}/n)^{1/2} }\big)$ has the same asymptotic distribution as $\sum_{i=1}^{n} (\lambda_i/n)^{1/2} \xi_{i}$.

Next we show that $n^{1/2} \text{tr}\big(({H_X-\frac 1n {\it\Omega}} )\tilde{S}^{1/2}\big)$ converges to 0 in probability.
Write
$$n^{1/2} \text{tr}\big(({H_X-\frac 1n {\it\Omega}} ) \tilde{S}^{1/2}\big) =\text{tr}(A\tilde{B}) =\sum_{i=1}^{m} \sum_{j=1}^{m} a_{ij} \tilde{b}_{ij},$$
where  $A=(a_{ij})_{m \times m}= (\frac{1}{n}X^\top X)^{-1}-\tilde\Delta^{-1} \stackrel{p}\longrightarrow 0 $ and  $\tilde{B}=(\tilde{b}_{ij})_{m \times m}=X^\top(\tilde{S}/n)^{1/2}X$. Denote $B_0=\tilde\Delta^{-1/2} \tilde{B} \tilde\Delta^{-1/2}$ and by Lemma \ref{lem3}
\begin{align*}
\text{tr}(B_0^\top B_0)=\text{tr}\big(\tilde\Delta^{-1} \tilde{B} \tilde\Delta^{-1} \tilde{B} \big)=\text{tr}\big({\it\Omega}(\tilde{S}/n)^{1/2}{\it\Omega}(\tilde{S}/n)^{1/2}  \big) \leq \Big(\text{tr}\big({\it\Omega}(\tilde{S}/n)^{1/2} \big)\Big)^{1/2},
\end{align*}
where $\text{tr}\big({{\it\Omega} (\tilde{S}/n)^{1/2} }\big)$  tends to $\sum_{i=1}^{\infty}  (\lambda^*_i)^{1/2}  \xi_i$ who is bounded in probability. It implies that $\text{tr}(B_0^\top B_0)$, i.e., the $(i,j)$th entry of $B_0$ is bounded in probability such that the $(i,j)$th entry of $ \tilde{B} =\tilde\Delta^{1/2} B_0 \tilde\Delta^{1/2}$ is as well bounded, $i,j=1,\cdots,m$. One can derive that  $n^{1/2} \text{tr}\big(({H_X-\frac 1n {\it\Omega}} )\tilde{S}^{1/2}\big)$  converges to 0 in probability. It follows that
$\text{tr}\big(H_X (n\tilde{S})^{1/2} \big)$ has the same asymptotic distribution as $\sum_{i=1}^{n}   (\lambda_i/n)^{1/2}  \xi_{i}$ as $n\rightarrow\infty$.

By the extension of Lemma \ref{lem2}, we have $\text{tr}\big(({\bf I}_n - H_X)\tilde{S}^{1/2}\big)/\sum_{i=1}^{n}\lambda_i^{1/2} \stackrel{p}\longrightarrow 1$ as $n\rightarrow\infty$. Consequently, $T_{sqrt}$ has the same asymptotic distribution with $m^{-1}\sum_{i=1}^{n} \eta_i \xi_{i}$.

%%%%%%%%  proof of Theorem 3 %%%%%%%%%

\subsection*{Proof of Theorem \ref{thm3}.}
For the linear kernel $\psi(\cdot,\cdot)$, it follows that
\begin{align*}
&\frac{1}{n}X^\top\tilde{S}X= \frac{1}{n} \sum_{i=1}^{n} \sum_{j=1}^{n} (x_i-\bar x)s_{ij}(x_j-\bar x)^\top\\
=&\Big[(n)^{-1}\sum_{i=1}^{n}  (x_i-\bar x) ( x_i^\top \bm \beta +  \varepsilon_i^\top)\Big] \Big[ \sum_{j=1}^{n}(\bm \beta^\top  x_j  +  \varepsilon_j)  (x_j-\bar x)^\top\Big].
\end{align*}
Let $G=(n)^{-1/2} \sum_{i=1}^{n}  (x_i-\bar x) ( x_i^\top \bm \beta+\varepsilon_i^\top)=\frac{1}{n^{1-\iota}} \sum_{i=1}^{n}  (x_i-\bar x)  x_i^\top \tilde{\boldsymbol{\beta}}_n+(n)^{-1/2} \sum_{i=1}^{n}  (x_i-\bar x) \varepsilon_i^\top\triangleq G_1+G_2$, where $G_2$ is bounded in probability by the Chebyshev's inequality and then $G_2/n^\iota$ are converge in probability to zero.  By the law of large numbers, it can be obtained that
$$Q_1=\frac 1n \sum_{i=1}^{n} s_{ij} = \frac{1}{n} \sum_{i=1}^{n}  ( x_i^\top \bm \beta +  \varepsilon_i^\top) (\bm \beta^\top  x_i  +  \varepsilon_i) \stackrel{p} \longrightarrow \Delta_\epsilon ~~\text{and}~~ \ Q_2=\frac{1}{n^2} \sum_{i=1}^{n} \sum_{j=1}^{n} s_{ij} \stackrel{p} \longrightarrow 0$$ and
the $(i,j)$th entry of $G_1/n^\iota$ converges in probability to the $(i,j)$th entry of $\Delta\tilde{\boldsymbol{\beta}}$ for $i=1,\cdots,m$ and $j=1,\cdots,q$. Then we have the numerator of $T$,
$$\frac{1}{n^{2\iota}}\text{tr}\big(H_X\tilde{S}\big)
=\frac{1}{n^{2\iota}}\text{tr}\big( (\frac{1}{n}X^\top X)^{-1} GG^\top\big)
\stackrel{p} \longrightarrow \text{tr}(\tilde\Delta^{-1} \Delta \tilde{\boldsymbol{\beta}}\tilde{\boldsymbol{\beta}}^\top \Delta^\top ) > 0.$$
and the denominator of $T_{pseudo}$,
$\text{tr}\big(\tilde{S}/n\big)=\text{tr}(Q_1-Q_2) \stackrel{p} \longrightarrow \text{tr}(\Delta_{\varepsilon})>0.$
It follows that
$$\frac{T_{pseudo}}{n^{2\iota}}=\frac{n-m}{mn}\frac{\text{tr}\big(H_X\tilde{S}/n^{2\iota} \big)}{\text{tr}\big(({\bf I}_n - H_X)\tilde{S}/n\big)}
\stackrel{p} \longrightarrow \frac 1m \frac{\text{tr}(\tilde\Delta^{-1} \Delta \tilde{\boldsymbol{\beta}} \tilde{\boldsymbol{\beta}}^\top \Delta^\top ) }{\text{tr}(\Delta_{\varepsilon})}>0$$
and then for any $\tilde\tau_0>0$,
$$\lim_{n\rightarrow\infty}  P\Big(\Big|\frac{T_{pseudo}}{n^{2\iota}}- \frac 1m \frac{\text{tr}(\tilde\Delta^{-1} \Delta \tilde{\boldsymbol{\beta}} \tilde{\boldsymbol{\beta}}^\top \Delta^\top ) }{\text{tr}(\Delta_{\varepsilon})}\Big|>\tilde\tau_0\Big)=0 $$

Now we show that
$P\big(T_{pseudo}> F_1^{-1}(1-\alpha)\big)=P\big(T_{pseudo}/n^{2\iota} > F_1^{-1}(1-\alpha)/n^{2\iota}\big)\stackrel{p} \longrightarrow 1$ as $n\rightarrow\infty$.

%%%%%%%%  proof of Theorem \ref{lem4} %%%%%%%%%

\subsection*{Proof of Theorem \ref{thm4}.}
Denote $A_0= H_X(\tilde{S}/n)^{1/2} $ whose eigenvalues are nonnegative by Lemma \ref{lem3} and then $\text{tr}\big(H_X(\tilde{S}/n)^{1/2})=\text{tr} (A_0) \geq  \text{tr} (A_0^\top A_0)^{1/2} =\text{tr}\big(H_X\tilde{S}/n) $ where $\text{tr}\big(H_X\tilde{S}/n)$ converges to a positive constant in probability by Theorem \ref{thm3}. It follows that
$$\frac{T_{sqrt}}{n}=\frac{\text{tr}\big(H_X\tilde{S}^{1/2}/mn \big)}{\text{tr}\big(({\bf I}_n - H_X)\tilde{S}^{1/2}/(n-m)\big)} \geq  \frac{n-m}{n}  \frac{\text{tr} (A_0^\top A_0)^{1/2} }{\sum_{i=1}^{n}(\lambda_i/n)^{1/2} }>0.  $$
Consequently, we have $P\big(T_{sqrt}> F_2^{-1}(1-\alpha)\big)=P\big(T_{sqrt}/n > F_2^{-1}(1-\alpha)/n\big)\stackrel{p} \longrightarrow 1$ as $n\rightarrow\infty$.

\end{document}